\documentclass{autart2}    

\usepackage{graphicx}      
\usepackage{natbib}        
\usepackage{graphics} 
\usepackage{times} 
\usepackage{amsmath} 
\usepackage{amssymb}  

\usepackage{tabu}
\usepackage{enumitem}
\usepackage{url}
\usepackage{mdframed}
\usepackage{algorithmic}
\usepackage{algorithm}

\def\QED{~\rule[-1pt]{5pt}{5pt}\par\medskip}
\newtheorem{postulate}{\bf Postulate}
\newenvironment{upd}{}

\begin{document}

\begin{frontmatter}

\title{Directed Information as Privacy Measure in Cloud-based Control\thanksref{footnoteinfo}} 

\thanks[footnoteinfo]{A preliminary version of this paper will be presented at ACC2017 \cite{tanaka2017directed}. }

%
\author[]{Takashi Tanaka}
\hspace{3ex}
\author[]{Mikael Skoglund}
\hspace{3ex}
\author[]{Henrik Sandberg}
 \hspace{3ex}
\author[]{Karl Henrik Johansson}

\address{School of Electrical Engineering, KTH Royal Institute of Technology, Stockholm, Sweden.}


\begin{abstract} 
  {\upd We consider cloud-based control scenarios in which clients
    with local control tasks outsource their computational or physical
    duties to a cloud service provider.  In order to address privacy
    concerns in such a control architecture, we first investigate the
    issue of finding an appropriate privacy measure for clients who
    desire to keep local state information as private as possible
    during the control operation. Specifically, we justify the use of
    Kramer's notion of causally conditioned directed information as a
    measure of privacy loss based on an axiomatic argument.
    Then we propose a methodology to design an optimal ``privacy
    filter'' that minimizes privacy loss while a given level of
    control performance is guaranteed. We show in particular that the optimal
    privacy filter for cloud-based Linear-Quadratic-Gaussian (LQG)
    control can be synthesized by a Linear-Matrix-Inequality (LMI)
    algorithm.  
The trade-off in the design is illustrated by a numerical example.}
\end{abstract}

\end{frontmatter}

\section{Introduction}

\begin{figure}[t]
    \centering
    \includegraphics[width=0.7\columnwidth]{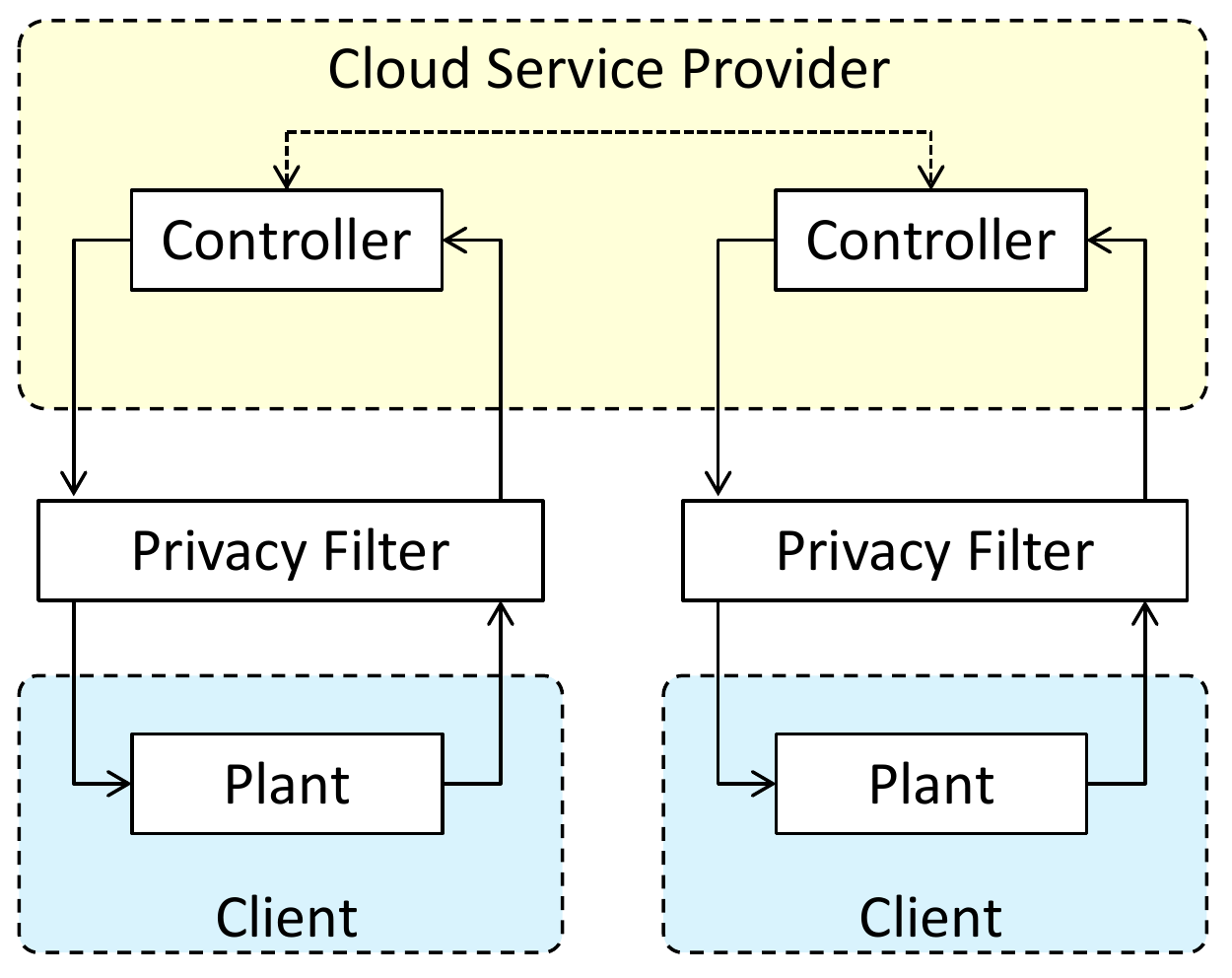}
    \caption{Cloud-based control.}
    \label{fig:cloud}
\end{figure}

{\upd
Leveraged by cloud computing technologies, the concept of cloud-based control has attracted much attention from industry in recent years.
Unlike the conventional situation where 
local agents are solely responsible for local control tasks, cloud-based control offers a flexible architecture in which a third party (i.e., cloud operator) provides control services (Fig.~\ref{fig:cloud}).
Advantages of such an architecture include the following.
\begin{enumerate}[leftmargin=5ex, label=(\roman*)] 
\item \label{cnt01} The local agent can outsource computational tasks to the cloud computer.
\item \label{cnt03} Global/shared information available to the cloud can improve control performances.
\item \label{cnt02} Physical resources needed for control actions can be provided by the cloud operator.
\item \label{cnt04} New kinds of services become available using data-mining technologies on large-scale operational data.
\end{enumerate}

Examples of cloud-based control strategies in category \ref{cnt01} include Model Predictive Control (MPC) of highly complex plants, where solving large-scale optimization problems in real-time is a critical requirement. For instance, \cite{hegazy2015industrial} studies an MPC-based operation of a large scale solar power plant, where the benefits of outsourcing computational tasks are discussed. 

A traffic monitoring and management system (e.g., \cite{hoh2012enhancing}) is an example of  cloud-based control in category \ref{cnt03}.
In this scenario, individual vehicles can be considered as clients, whose control tasks are to arrive at the destinations efficiently.
Since vehicles share a common infrastructure, the overall control performance is drastically improved by the existence of a centralized decision coordinator.

Cloud-based control services in category \ref{cnt02} provide not only
computational services but also physical resources. For instance, the
shared Energy Storage Systems (ESS) for smart grids
\cite{wang2013active,rahbar2016shared} can be considered as a
cloud-based control systems in this category.  In this example, clients (e.g.,
individual households) with unreliable renewable energy sources store
their excess energy in the shared ESS (operated by the cloud). While
such an architecture is reported to have cost advantages over
distributed ESS systems \cite{rahbar2016shared}, this introduces new
privacy risks since individual power consumption profiles can be
observed by the ESS operator.  As we will see, this is an important
example of cloud-based control in which privacy cannot be fully
protected by data encryption due to the actual physical signals involved (e.g.,
power consumption).

Finally, category \ref{cnt04} includes the concept of predictive
  manufacturing \cite{lee2013recent}.  This is an idea of collecting
and analyzing large amounts of operational data from machines in
production lines, targeting at improving productivity and safety by
predicting failures before they occur.  }

\subsection{Privacy concerns in cloud-based control}

{\upd While the concept of cloud-based control enhances conventional
  control technologies in various ways, it also brings new risks that
  did not exist in traditional scenarios.  Clearly, careless
  installation of cloud-based control strategies endangers the
  clients' privacy, since the architecture allows the cloud operator
  to learn sensitive information belonging to the clients, from the
  operational records.  Since private information can be a valuable asset 
  in our modern society, a reasonable assumption is that the cloud
  operator is intrinsically inclined to do so. Hence in this paper, we
  consider the cloud operator as a semi-honest (i.e., ``honest but
  curious'') agent, meaning that it persistently tries to infer the
  clients' private information while executing the designated control
  algorithm faithfully.

  Protecting privacy in cloud-based control requires multiple layers
  of data security technologies, such as data encryption and data
  perturbation.  The latter technology has been actively studied in
  the database literature in recent years, where the trade-off between
  data utility and privacy with respect to various metrics (e.g.,
  differential privacy \cite{dwork2008differential}, $k$-anonymity
  \cite{sweeney2002k}, information theoretic privacy
  \cite{du2012privacy}) is thoroughly studied.  However,
  these privacy metrics were introduced in control theoretic settings
  relatively recently (\cite{cortes2016differential} and references
  therein), and it is safe to say that privacy concerns in cloud-based
  control is an area still in it infancy.

  There are several reasons why the existing privacy mechanisms cannot
  be (and should not be) naively used in cloud-based control.  First,
  data encryption technologies alone (e.g., full or partial
  homomorphic encryption
  \cite{kogiso2015cyber,farokhi2016secure,shoukry2016privacy}, data
  obfuscation \cite{wang2011secure}, multi-party computation schemes
  \cite{du2001secure}) may not be sufficient to protect privacy.
  Their limitations come from several reasons, including (i) as of
  today the computational requirement for encryption is still far from
  practical (e.g., homomorphic encryption
  \cite{farokhi2016secure,shoukry2016privacy}), (ii) some encryption
  technologies need public keys to be delivered reliably, but this
  itself requires separate security guarantees, and (iii) in some
  situations the cloud operator inevitably has access to
  decrypted/physically meaningful data.  To see the last point, recall
  the aforementioned shared ESS example, where the power inflow to
  an individual household (which contains sensitive information) is
  inevitably observable by the ESS operator.  While data encryption
  cannot be used here, notice that data perturbation can be applied to
  physical signals as well. For instance, battery load hiding
  (using, e.g., local energy storage devices \cite{tan2013increasing, chin2016privacy})
   is a data perturbation technique to enhance
  smart meter privacy.

  Second, privacy notions targeting at single-stage data disclosure
  mechanisms (which are often the case in the database literature) are
  in general not sufficient to accommodate privacy issues in
  cloud-based control. A good privacy notion must  respect the
  fact that privacy leakage occurs over multiple time steps, and the
  data from the past stored in the cloud can potentially be used to
  threaten privacy at the present time.  Also, the existence of
  information feedback must be carefully taken into account.  Namely,
  the cloud has certain influences (through control inputs) on the
  future private information (state of client plants), and hence
  appropriate statistical conditioning is needed to distinguish
  private information from public information.  }

Finally, in general, establishing the adequacy of privacy notions
(e.g., differential privacy, $k$-anonymity, information theoretic
privacy) in a particular application (in our case, cloud-based control)
often requires subtle examination of the context.  In fact, many of
the available privacy notions and their validity are sensitive,
explicitly or implicitly, to the problem setting at hand and premises where
those notions were originally introduced. For instance,
\cite{machanavajjhala2007diversity} shows by simple counterexamples
that $k$-anonymity is fragile against side-information.  {\upd While
  differential privacy is shown to be stronger than
  information-theoretic privacy in a certain sense \cite{de2012lower},
  it is also demonstrated in a somewhat different scenario that it
  does not provide any guarantee in information-theoretic privacy
  \cite{du2012privacy}.

  To cope with these difficulties, and to embellish privacy
  discussions for cloud-based control, this paper introduces an
  axiomatic approach to identify an appropriate privacy notion for
  cloud-based control.}  Specifically, we propose a set of postulates,
which is a set of natural properties to be satisfied by a reasonable
notion of privacy in cloud-based control, and show that a particular
function, namely Kramer's notion of \emph{causally conditioned
  directed information}, arises as a unique candidate.  An axiomatic
characterization also provides a convenient interface between the theory
and practice of privacy considerations.  As discussed above, it is often
difficult to judge whether a \emph{given} notion of privacy is
appropriate for individual applications.  In contrast, axioms are
often easier to discuss in practical contexts.  Axioms also provide a
solid mathematical basis on which rigorous theory of privacy can be
developed.  {\upd As a privacy protection mechanism, we propose an
  additional layer (\emph{privacy filter}) bridging the cloud and
  clients, which has a dedicated role to control the leakage of
  private data (Fig.~\ref{fig:cloud}).  We propose joint design of
  privacy filter and control algorithms, so that the overall system is
  able to balance utility of cloud-based control and privacy losses
  (with respect to the derived privacy notion). }

\subsection{Related work}

Privacy has been  extensively studied in the database literature in recent years.
While ad hoc approaches for privacy (sub-sampling, aggregation, and suppression) have a long history, one of the first formal definitions of privacy is given by $k$-anonymity \cite{sweeney2002k}. Extensions of this notion include $t$-closeness and $l$-diversity \cite{machanavajjhala2007diversity}. Differential privacy \cite{dwork2008differential} has been particularly popular since its introduction, partly because of its convenient property that no prior on the database content is needed nor used. 
Information-theoretic privacy in database is  considered in \cite{du2012privacy} and \cite{Poor}.

Privacy has only relatively recently become a topic of concern in the control-engineering literature. Some of the first works in the area treated consensus algorithms, and how participating agents can maintain some level of privacy despite sharing information with neighbors, see \cite{huang+12,manitara+13,mo+16}. Differential privacy, which 
 was originally developed for database privacy, can quite generally be adapted to a control-theoretic context as shown in \cite{ny+14,wang+14}, and also in particular filtering and control applications, e.g.,  \cite{huang+12,sandberg+15,mo+16,nozari+16}. Based on game theory, alternative rigorous notions of privacy in a control and filtering context have been obtained in \cite{akyol+15,farokhi+15}.

A general introduction to information-theoretic security, secrecy, and
privacy can be found in \cite{BlochBarr}. Information theory has been
used to analyze various aspects of privacy in several different
problems settings. The problem of private information-retrieval
\cite{Sudan} was considered for example in \cite{Stadler} (and
references therein). Recent work reported in \cite{Jafar}
introduced the notion of capacity of private information-retrieval,
and characterized corresponding fundamental bounds.
Information-theoretic tools have also been utilized in the context of
differential privacy \cite{Andres,Kopf}. Very recent work summarized in
\cite{Zhang} studies the relation between differential privacy and
privacy quantified in terms of mutual information. This paper also
relates these two notions to the concept of identifiability.  The work
reported in \cite{Poor} introduced a general framework for
establishing a relation between privacy and utility based on
rate--distortion arguments. Similarly, \cite{Willems-bio} developed
analytic tools to support the characterization of leakage of privacy
in biometric systems.

\subsection{Contribution of this paper}

Contributions of this paper are summarized as follows.
\begin{itemize}[leftmargin=2ex]
\item[(a)] We provide a set of postulates (Postulates~\ref{post1}-\ref{post4}) characterizing basic properties of a privacy measure in cloud-based control, and show that Kramer's causally conditioned directed information arises as a unique candidate.
\item[(b)] We formulate an optimization problem characterizing optimal joint control and privacy filter policies, and derive its explicit solution in the LQG case. 
\end{itemize}

{\upd
We note that contribution (a) is crucially dependent on the recent result \cite{jiao2015justification} where a justification of the logarithmic loss functions is given via the so-called data-processing axiom.
(Note that such an axiomatic characterization of information measures has long and rich history  \cite{csiszar2008axiomatic}.)
Our contribution is an extension of \cite{jiao2015justification} to the multi-stage data disclosure mechanisms with information feedback, and its re-interpretation as a privacy axiom.

Notice also that, although we show that the causally conditioned directed information is the only candidate satisfying the considered set of postulates, we do \emph{not} claim that the considered postulates are the only possible characterization of privacy.
In fact, it is our important future work to examine carefully, possibly using real-world incidents of privacy attacks, whether the considered privacy postulates are appropriate or not. At the same time it is worth studying how a different set of postulates leads to a different notion of privacy.

We also note that axiomatic consistency is not the only criterion that
determines usefulness of various privacy notions.  For instance, to
design a privacy filter according to our privacy notion we need to
have precise knowledge of the system model (e.g., distributions of
process noises). This is a weakness compared to  mechanisms based
on differential privacy, which do not require prior knowledge of the
system.\footnote{On the other hand, in control we often have some prior knowledge of the system, which should be incorporated in the privacy filter design.}  }

\subsection{Notation}
Random variables are indicated by upper case symbols such as $X$. We denote by $P_X$, $P_{X,Y}$ and $P_{X|Y}$ the probability distribution of $X$,  the joint probability distribution of $X$ and $Y$, and the stochastic kernel of $X$ given $Y$, respectively. We use notation $P_{X|y}$ to emphasize that it is the conditional probability distribution of $X$ given $Y=y$. We write $H(X|z)$ and $I(X;Y|z)$ to denote the entropy and mutual information evaluated under $P_{X,Y|z}$, and define conditional entropy and conditional mutual information by $H(X|Z):= \mathbb{E}_{P_Z} H(X|z)$ and $I(X;Y|Z):= \mathbb{E}_{P_Z} I(X;Y|z)$. If $f$ is a function of a random variable $X$, denote by $\mathbb{E}_{P_X} f(X)$ or $\mathbb{E}_{P_X} f(x)$ the mathematical expectation. The cardinality of a set $\mathcal{X}$ is denoted by $|\mathcal{X}|$. A positive definite (resp. semidefinite) matrix $M$ is indicated by $M\succ 0$ (resp. $M \succeq 0$).

\section{Problem setting}
In this paper, a cloud-based control system is modeled by a discrete-time nonlinear stochastic control system.
We say that a random variable is \emph{public at time} $t$ if its realization is known to the cloud operator at time $t$. 
In contrast, by \emph{private random variable at time} $t$, we refer to random variables that the client wishes to keep confidential (in an appropriate sense discussed below) at time $t$.\footnote{According to this definition, note that random variables are public or private, or neither.} This classification reflects our premise that the cloud operator is semi-honest.
In this paper, we treat the state sequence $X^t\triangleq (X_1, ... , X_t)$ of the local plant up to time $t$ as the private random variable at time $t$.
We wish to introduce an appropriate measure of privacy loss that occurs during the operation of cloud-based control over a period $1\leq t\leq T$.

\begin{figure}[t]
    \centering
    \includegraphics[width=0.95\columnwidth]{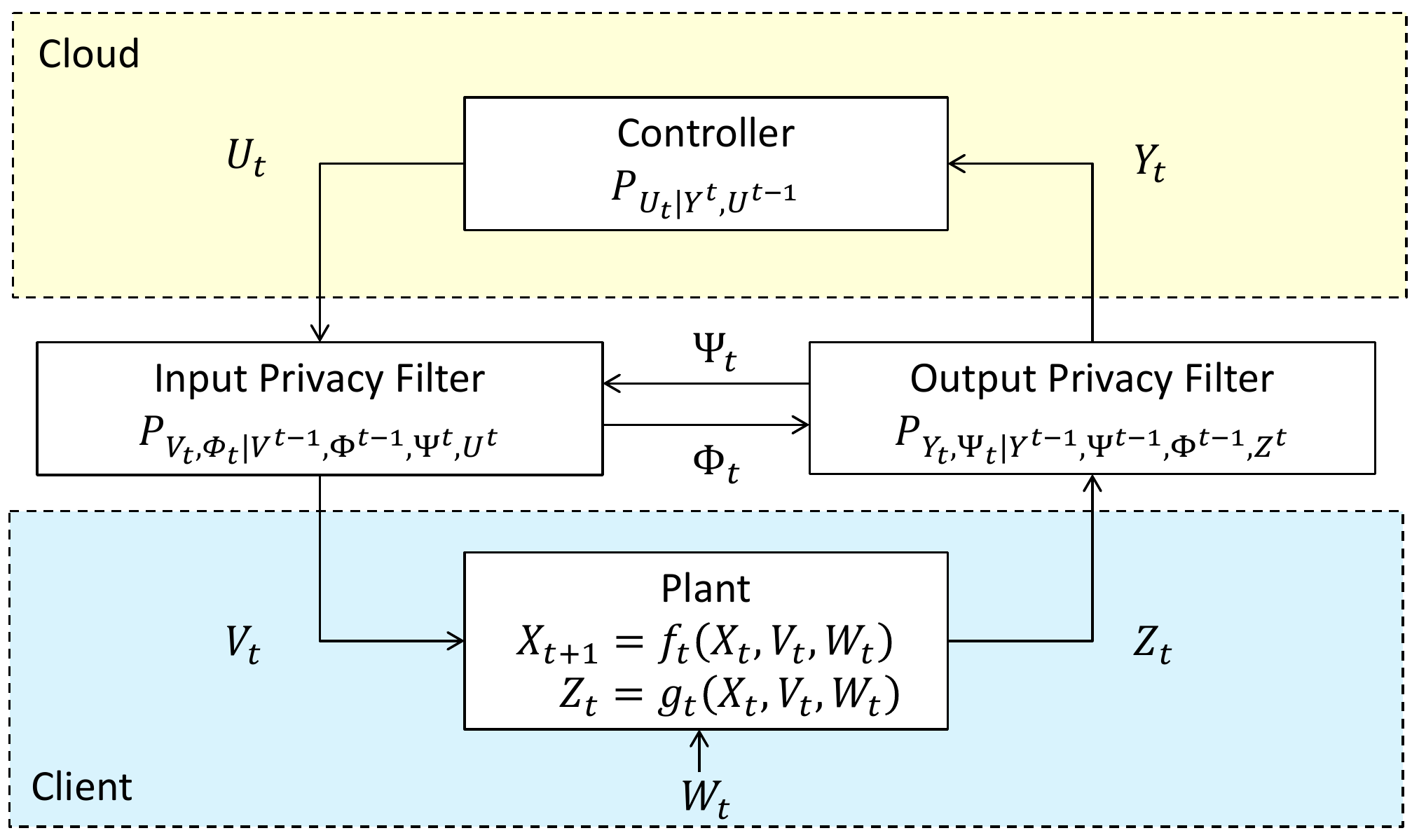}
    \caption{Privacy filter: General model.}
    \label{fig:inout}
\end{figure}

\begin{figure}[t]
    \centering
    \includegraphics[width=0.95\columnwidth]{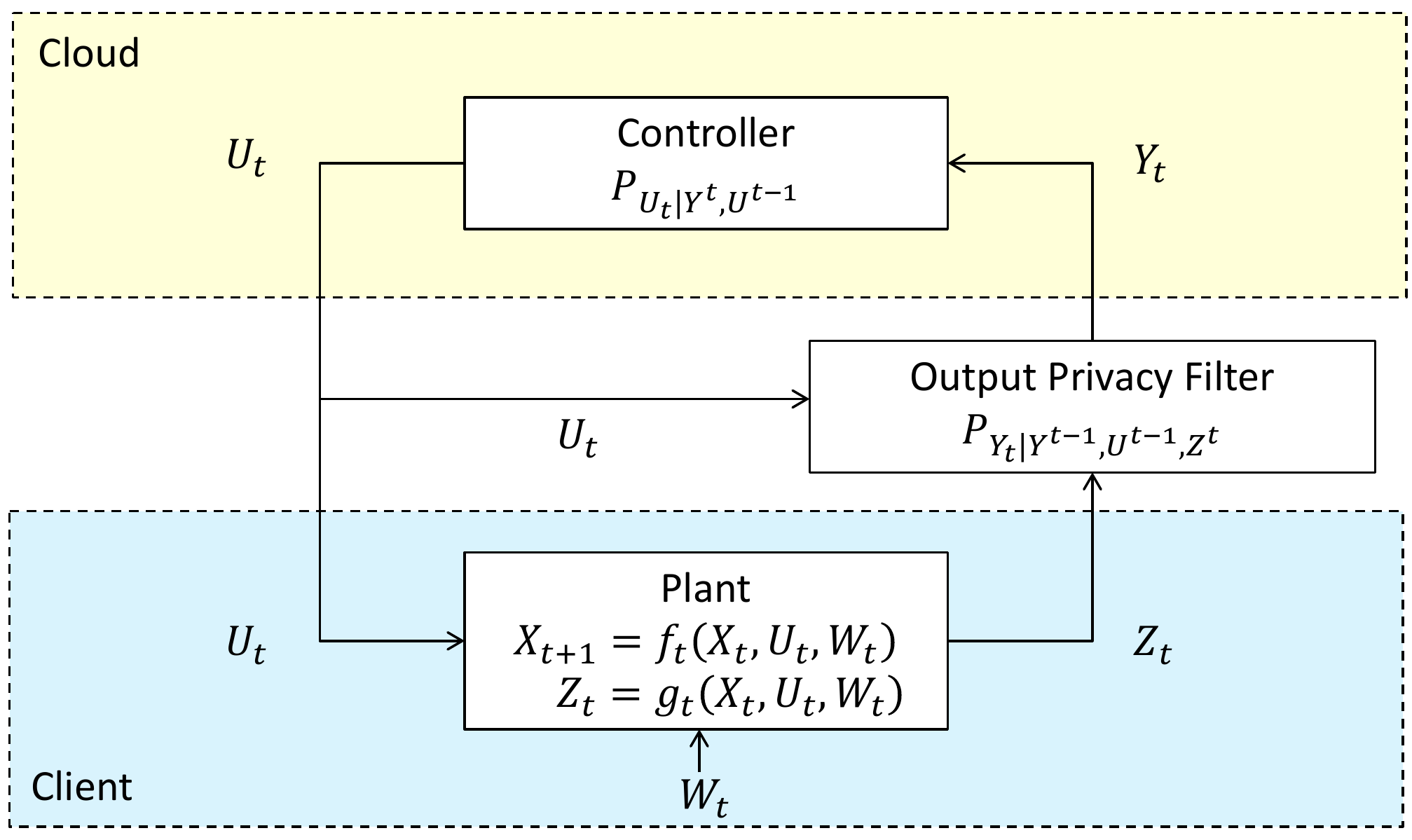}
    \caption{Privacy filter: Output filter only.}
    \label{fig:out}
\end{figure}

Fig.~\ref{fig:inout} illustrates the general structure of the class of
privacy filters considered.  An output filter prevents raw sensor data
to be disclosed to the cloud. An input privacy filter replaces the
control input $U_t$ with a different value $V_t$ to enhance privacy.
In general, the input and output filters can communicate with each
other via messages $\Psi_t$ and $\Phi_t$.  Privacy filters and
controller algorithms are in general randomized policies and have
memories of the past observations.  Thus, we model them as stochastic
kernels of the forms specified in Fig.~\ref{fig:inout}.
Fig.~\ref{fig:out} shows a simpler form of a privacy filter in which
the control input commanded by the cloud is directly applied to the
plant.  Since there is no input filter, this architecture is easier to
implement. For the rest of the paper, we focus on this simple
architecture in Fig.~\ref{fig:out}, and discuss privacy notions and
privacy filter design problems exclusively for this architecture.  In
Section~\ref{secaxiom}, we characterize our privacy notion
axiomatically, and then formulate a joint controller and output
privacy filter design problem in Section~\ref{secdesign}.  We derive
an optimal form of joint controller and output privacy filter in the
LQG regime in Section~\ref{seclqg}.

\section{Axiomatic characterization of privacy}
\label{secaxiom}
A meaningful notion of privacy must satisfy some basic properties.  In
this section, we first consider a single-stage data disclosure
mechanism and show that the only candidate function that satisfies the
natural set of postulates (axioms) is Shannon's mutual information
between private and published random variables.  Our arguments are
aligned with the development in \cite{jiao2015justification}, where
 mutual information arises as a unique function that characterizes
the value of side information in inference problems.  The set of
axioms used there is simple, and thus we argue that it can be
naturally used as a set of axioms for privacy.  Then, we apply this
observation to multi-stage feedback control systems and show the
unique candidate characterizing privacy loss in cloud-based control in
a satisfactory manner is the \emph{causally conditioned directed
  information}.

\subsection{Single-stage case}
Suppose $X$ and $Y$ are $\mathcal{X}$- and $\mathcal{Y}$-valued random variables with joint distribution $P_{X,Y}$. We temporarily assume that $\mathcal{X}$ and $\mathcal{Y}$ are countable sets, and denote by $\mathcal{P}_\mathcal{X}$ the space of probability distributions on $\mathcal{X}$.
Assuming that $X$ is a private random variable, we wish to quantify the privacy loss due to the disclosure of a random variable $Y$.

First, we quantify the ``hardness" of inferring $X$ using the notion of \emph{loss function}. Generally speaking, a random variable $X$ is hard to infer if the expected posterior value of observation (i.e., the degree of ``surprise" that occurs when observing a realization $x\in \mathcal{X}$) cannot be made small.
The posterior value of the observation is a function of the observed realization $x\in \mathcal{X}$ and a prior distribution $Q_X\in\mathcal{P}_{\mathcal{X}}$ assumed by the observer. We refer to such a function $\ell: \mathcal{X}\times \mathcal{P}_{\mathcal{X}}\rightarrow \mathbb{R}$ as a \emph{loss function}. In the literature, it is also called the scoring rule \cite{gneiting2007strictly} or self-information \cite{merhav1998universal}.
Note that for a given choice of $\ell$, the task of inference is to minimize $\ell$ by properly assuming $Q_X$.
Among many options, the \emph{logarithmic loss function} $\ell(x, Q_X)=\log \frac{1}{Q_X(x)}$ is frequently used in the literature. We will motivate this choice later in this paper (Postulate~\ref{post2}).

Let  $P_X \in \mathcal{P}_{\mathcal{X}}$ be the true probability distribution of $X$, and $Q_X \in \mathcal{P}_{\mathcal{X}}$ be the assumed distribution. In general $P_X\neq Q_X$.
If $S(P_X, Q_X) \triangleq \mathbb{E}_{P_X} \ell(X, Q_X)$, the quantity $\inf_{Q_X} S(P_X, Q_X)$ is referred to as the \emph{Bayes envelope}.
A loss function is said to be proper if 
$
\inf_{Q_X} S(P_X, Q_X)=S(P_X, P_X)$.
It can easily be shown that the logarithmic loss function is proper, and the associated Bayes envelope coincides with the entropy $H(X)$ of $X$:
\[
\inf_{Q_X} S(P_X, Q_X)=S(P_X, P_X)=\mathbb{E}_{P_X} \log\frac{1}{P_X(x)}=H(X).
\]

Now we introduce the first postulate characterizing our privacy notion. It states that privacy loss due to disclosing $Y$ is measured
by the expected difference in Bayes envelope evaluated before and after observing $Y$.
Given a loss function $\ell$ and the joint distribution $P_{X,Y}$, we refer to the privacy loss evaluated this way as the \emph{privacy leakage function}, and denote it by $L(\ell, P_{X,Y})$.

\begin{postulate}
\label{post1}
The privacy leakage function is in the form
\begin{align*}
&L(\ell,P_{X,Y})= \\
& \;\; \inf_{Q_X} \mathbb{E}_{P_{X}} [\ell(X, Q_X)]-\mathbb{E}_{P_Y}\inf_{Q_{X|Y}} \mathbb{E}_{P_{X|Y}} [\ell(X, Q_{X|Y})].
\end{align*}
\end{postulate}
The first term on the right hand side is the Bayes envelope evaluated without side information $Y$, while in the second term, the assumed distribution $Q_{X|Y}$ is allowed to depend on $Y$.
Hence, $L(\ell,P_{X,Y})$ is understood to be the improvement in the estimation quality due to the side information $Y$.
If the loss function $\ell$ is logarithmic, the privacy leakage function $L(\ell,P_{X,Y})$ defined above coincides with the mutual information between $X$ and $Y$, i.e.,
\[ 
L(\ell, P_{X,Y})=H(X)-H(X|Y)=I(X;Y).
\]
Up to now, the logarithmic loss function is just an example among many other possible choices of loss functions. It turns out that it is the only option that satisfies the following natural postulate.

\begin{postulate}\label{post2}(Data-processing axiom \cite{jiao2015justification})
For any distribution $P_{X,Y}$ on $\mathcal{X}\times \mathcal{Y}$, the information leakage function $L(\ell,P_{X,Y})$ satisfies
\begin{equation}
\label{dpaxiom}
L(\ell,P_{T(X),Y}) = L(\ell,P_{X,Y})
\end{equation}
for every $T: \mathcal{X}\rightarrow \mathcal{X}$ such that $T(X)$ is a sufficient statistic of $X$ for $Y$, i.e., the following Markov chains\footnote{$X$ -- $Y$ -- $Z$ means that $X$ and $Z$ are conditionally independent given $Y$.} hold:
\begin{equation}
\label{post2mc}
T(X) \text{--} X \text{--}  Y, \;\;\; X \text{--} T(X) \text{--} Y.\end{equation}
In \eqref{dpaxiom}, the joint distribution $P_{T(X),Y}$ on $\mathcal{X}\times \mathcal{Y}$ is defined by
\[ P_{T(X),Y}(T(\mathcal{B}_\mathcal{X})\times \mathcal{B}_\mathcal{Y})=P_{X,Y}(\mathcal{B}_\mathcal{X}\times \mathcal{B}_\mathcal{Y})\]
for all subsets $\mathcal{B}_\mathcal{X}$ and $\mathcal{B}_\mathcal{Y}$ of $\mathcal{X}$ and $\mathcal{Y}$, respectively.
\end{postulate}

{\upd
\begin{rem}
For instance, one can think of $\mathcal{X}=\mathbb{R}^n$ and $T$ being a coordinate transformation of $X$. The identity \eqref{dpaxiom} then implies that the privacy loss is uniquely defined no matter what coordinate is chosen to represent private random variables. Postulate~\ref{post2} argues that it would be natural to require \eqref{dpaxiom} for any transformation satisfying \eqref{post2mc}.
\end{rem}
In \cite{jiao2015justification}, a weaker axiom with an inequality
\begin{equation}
\label{dpaxiomineq}
L(\ell,P_{T(X),Y}) \leq L(\ell,P_{X,Y})
\end{equation}
is used. However, for our purpose \eqref{dpaxiom} and
\eqref{dpaxiomineq} are equivalent, and have the same consequences.
Although the requirement of Postulate~\ref{post2} seems rather mild,
it has strong implications as summarized by the next theorem.  }
\begin{thm}\label{theo1} (Justification of mutual information \cite{jiao2015justification})
Let $\mathcal{X}$ be a finite set with $|\mathcal{X}|\geq 3$.
Under Postulate~\ref{post2}, the privacy leakage function is uniquely determined by the mutual information
\[
L(\ell, P_{X,Y})=I(X;Y)
\]
up to a positive multiplicative factor.
\end{thm}
\begin{pf}
{\upd
A complete proof is given in \cite{jiao2015justification}. Notice that the steps in which the inequality version of the axiom \eqref{dpaxiomineq} is used (namely, equations (24) and (80) in \cite{jiao2015justification}) can also be established by the equality version of the axiom \eqref{dpaxiom}. \hfill \QED
}
\end{pf}
\begin{rem}
The result of Theorem~\ref{theo1} can be extended to the case with continuous random variables $X$ and $Y$ using the formula \cite[Ch. 2.5]{gallager1968information}, \cite[Ch. 3.5]{pinsker1964information},\cite[Ch. 7.1]{gray1990entropy}: 
\begin{equation}
\label{eqquantizedmi}
I(X;Y)=\sup I([X], [Y]).
\end{equation}
The right-hand-side of \eqref{eqquantizedmi} denotes the supremum of mutual information between discrete random variables $[X]$ and $[Y]$ over all finite quantizations. If we consider a supremum achieving sequence of quantizers, and require the data-processing axiom to be satisfied by each element of the sequence, we obtain $I(X;Y)$ as the unique privacy leakage function for  continuous random variables $X$ and $Y$.
\end{rem}

\subsection{Multi-stage case}
Based on the single-stage discussion in the previous section, in this section we propose a multi-stage privacy measure suitable for cloud-based control (Fig.~\ref{fig:out}). To proceed, we introduce the following additional postulates.
\begin{postulate}
\label{post3}
The private random variable at time $t$ is $X^t$, while $Y^{t-1}$ and $U^{t-1}$ are public at time $t$.
\end{postulate}

We first characterize the instantaneous privacy loss at time step $t$ due to the disclosure of $Y_t$. By Postulate~\ref{post3}, we need to characterize the privacy leakage function for $X^t$ due to disclosing $Y_t$ under the joint distribution $P_{X^t,Y_t|y^{t-1}, u^{t-1}}$. Notice that, by Postulate~\ref{post3}, $y^{t-1}$ and $u^{t-1}$ are public knowledge.

Let $\ell$ be a loss function as in the preceding subsection. For every realization $(y^{t-1}, u^{t-1})$, Postulate~\ref{post2} requires the privacy leakage function $L(\ell, P_{X^t,Y_t|y^{t-1}, u^{t-1}})$ to satisfy
\[
 L(\ell,P_{T(X^t), Y_t|y^{t-1}, u^{t-1}}) = L(\ell,P_{X^t, Y_t|y^{t-1}, u^{t-1}})
\]
whenever $T(X^t)\in \mathcal{X}^t$ is a sufficient statistic of $X^t$ for $Y_t$ given $Y^{t-1}=y^{t-1}$ and $U^{t-1}=u^{t-1}$, i.e., the following Markov chains hold under $P_{X^t, Y_t|y^{t-1}, u^{t-1}}$:
\[ T(X^t) \text{--} X^t \text{--} Y_t, \;\;\; X^t \text{--} T(X^t) \text{--} Y_t.\]
From Theorem~\ref{theo1}, we conclude that the only loss function (up to positive multiplicative factors) that satisfies the above equality is the logarithmic one, and with necessity we have 
\[ L(\ell, P_{X^t, Y_t|y^{t-1}, u^{t-1}}) =I(X^t; Y_t|y^{t-1}, u^{t-1}). \]
Thus, if $Y^{t-1}$ and $U^{t-1}$ have a joint distribution $P_{Y^{t-1},U^{t-1}}$, the expected privacy loss at time step $t$ is
\[ \mathbb{E}_{P_{Y^{t-1},U^{t-1}}}L(\ell, P_{X^t, Y_t|y^{t-1}, u^{t-1}})=I(X^t; Y_t|Y^{t-1}, U^{t-1}).
\]
Finally, we assume that our privacy notion satisfies the following natural property.
\begin{postulate}
\label{post4}
The expected total privacy loss over the horizon $t=1,2, ... , T$ has a stage-additive form over the expected instantaneous privacy losses.
\end{postulate}

Under Postulate~\ref{post4}, the expected total privacy loss is 
\begin{equation}
\label{eqconddi}
\sum_{t=1}^T I(X^t;Y_t|Y^{t-1}, U^{t-1})=: I(X^T\rightarrow Y^T \| U^{T-1}).
\end{equation}
The notation on the right hand side of \eqref{eqconddi} is introduced in \cite{kramer2003capacity}. We refer to this quantity as Kramer's causally conditioned directed information. Thus, we obtain:
\begin{prop}
Under Postulates~\ref{post1}-\ref{post4}, Kramer's causally conditioned directed information $I(X^T\rightarrow Y^T \| U^{T-1})$ is the only function (up to positive multiplicative factors) quantifying the expected privacy loss in the cloud-based control in Fig.~\ref{fig:out}.
\end{prop}

\section{Privacy-preserving cloud-based control design}
\label{secdesign}
Suppose that the performance of the cloud-based control system is
measured by a stage-wise additive cost function
$\sum_{t=1}^T \mathbb{E}c(X_{t+1},U_t)$.  Then, privacy loss in 
cloud-based control with a given control performance requirement
$\delta$ is minimized by solving
\begin{subequations}
\label{minIstC}
\begin{align}
\min \;\; & I(X^T\rightarrow Y^T \| U^{T-1}) \\
\text{s.t.} \;\;\; & \sum_{t=1}^T \mathbb{E} c(X_{t+1},U_t) \leq \delta.
\end{align}
\end{subequations}
Likewise, the best achievable control performance under the privacy constraint is characterized by flipping the constraint and objective functions in \eqref{minIstC}.
In both cases, the optimization domain is the space of the sequence of Borel measurable stochastic kernels
\begin{equation}
\label{eqoptdomain}
\mathcal{D}=\{P_{U_t|Y^t, U^{t-1}}\; , \;\;P_{Y_t|Y^{t-1}, U^{t-1}, Z^t} \}_{t=1}^T
\end{equation}
characterizing joint controller and output privacy filter policies.\footnote{We consider $P_{U_1|Y^1, U^{0}}=P_{U_1|Y^1}$ and $P_{Y_1|Y^{0}, U^{0}, Z^1}=P_{Y_1|Z^1}$.}
Since \eqref{minIstC} is an infinite dimensional optimization problem, it is in general  difficult to obtain an explicit form of an optimal solution. In Section~\ref{seclqg}, we consider a special case in which it is possible. 
{\upd For the later use, we next present a fundamental inequality showing that the privacy leakage is at least $I(X^T\rightarrow U^T):= \sum_{t=1}^T I(X^t;U_t|U^{t-1})$.
\begin{lem}
\label{lemfbineq}
For all joint controller and output privacy filter policies in \eqref{eqoptdomain}, we have
\[
I(X^T\rightarrow U^T) \leq I(X^T\rightarrow Y^T \| U^{T-1}).
\]
\end{lem}
\begin{pf}The inequality is directly verified as follows.
\begin{align*}
&I(X^T\rightarrow Y^T \| U^{T-1})-I(X^T\rightarrow U^T) \\
=&\sum\nolimits_{t=1}^T \left[I(X^t;Y_t|Y^{t-1},U^{t-1})-I(X^t;U_t|U^{t-1}) \right] \\
\stackrel{(a)}{=}&\sum\nolimits_{t=1}^T \left[I(X^t;Y_t, U_t|Y^{t-1},U^{t-1})-I(X^t;U_t|U^{t-1}) \right] \\
\stackrel{(b)}{=}&\sum\nolimits_{t=1}^T \left[I(X^t;Y^t|U^t)-I(X^t;Y^{t-1}|U^{t-1}) \right] \\
\stackrel{(c)}{=}&\sum\nolimits_{t=1}^T \left[I(X^t;Y^t|U^t)-I(X^{t-1};Y^{t-1}|U^{t-1}) \right] \\
\stackrel{(d)}{=}&I(X^T;Y^T|U^T) \geq 0.
\end{align*}
Equality (a) holds since $I(X^t;Y_t, U_t|Y^{t-1},U^{t-1})=I(X^t;Y_t|Y^{t-1},U^{t-1})+I(X^t;U_t|Y^t,U^{t-1})$, and the second term is zero since $X^t$ and $U_t$ are conditionally independent given $(Y^t,U^{t-1})$. To see (b), apply the chain rule for the mutual information in two different ways:
\begin{align*}
&I(X^t;Y^t,U_t|U^{t-1}) \\
&=I(X^t;Y^{t-1}|U^{t-1})+I(X^t;Y_t,U_t|Y^{t-1}, U^{t-1}) \\
&=I(X^t;U_t|U^{t-1})+I(X^t;Y^t|U^t).
\end{align*}
Equality (c) holds as $I(X^t;Y^{t-1})=I(X^{t-1};Y^{t-1}|U^{t-1})+I(X_t;Y^{t-1}|X^{t-1},U^{t-1})$, and the second term is zero since $X_t$ and $Y^{t-1}$ are conditionally independent given $(X^{t-1}, U^{t-1})$.
Finally, telescoping cancellations of terms show (d).
\hfill \QED
\end{pf}
}

So far we have provided a justification of $I(X^T\rightarrow Y^T \| U^{T-1})$ as a measure of privacy loss. Next, we discuss how this quantity imposes a fundamental limitation in estimating private random variables.

\subsection{Implication via distortion-rate function}
Consider an optimal joint controller and output privacy filter policy solving \eqref{minIstC}, and let $\gamma$ be the optimal value.
 By Postulate~\ref{post4}, the total privacy loss $\gamma$ can be written as $\gamma=\sum_{t=1}^T \gamma_t$, where
\begin{equation}
\label{eqlossatt}
\gamma_t=I(X^t;Y_t|Y^{t-1}, U^{t-1}) \;\; (\geq 0)
\end{equation}
is the privacy loss at time $t$. 
To see how \eqref{eqlossatt} guarantees privacy  against inferring $X^t$  at time $t$ even after disclosing $Y_t$, consider an estimate of $X^t$ of the form 
$
\hat{X}^t: \mathcal{Y}^t\times \mathcal{U}^{t-1}\rightarrow \mathcal{X}^t.
$
Since realizations of $Y^{t-1}$ and $U^{t-1}$ are prior knowledge at time $t$, $\hat{X}^t$ can be viewed as a function of $Y_t$ alone, and thus $X^t$ -- $Y_t$ -- $\hat{X}^t$ forms a Markov chain given $(Y^{t-1}, U^{t-1})$. By the data-processing inequality,
\[
I(X^t;\hat{X}^t|Y^{t-1}, U^{t-1}) \leq  I(X^t;Y_t|Y^{t-1}, U^{t-1})=\gamma_t.
\]
In other words, the expected mutual information between $X^t$ and $\hat{X}^t$ is bounded by $\gamma_t$: 
\begin{equation}
\label{eqgammat}
\mathbb{E}_{P_{Y^{t-1},U^{t-1}}} I(X^t;\hat{X}^t|y^{t-1}, u^{t-1}) \leq \gamma_t.
\end{equation}
This inequality imposes a fundamental limitation of estimation accuracy in the following sense.
 Let $\rho_t: \mathcal{X}^t \times \mathcal{X}^t \rightarrow [0,\infty)$ be an arbitrary distortion function. 
For a given source distribution $P_{X^t|y^{t-1}, u^{t-1}}$, let $D_t: [0,\infty)\rightarrow [0,\infty)$ be the distortion-rate function \cite{CoverThomas}.
By definition of the distortion-rate function, for any joint distribution $P_{X^t,\hat{X}^t|y^{t-1}, u^{t-1}}$, we have
\[
\mathbb{E}_{P_{X^t,\hat{X}^t|y^{t-1}, u^{t-1}}} \rho_t(X^t,\hat{X}^t) \geq D_t(I(X^t;\hat{X}^t|y^{t-1}, u^{t-1})).
\]
Taking expectation with respect to $P_{Y^{t-1}, U^{t-1}}$, we have
\begin{align*}
\mathbb{E}\rho_t(X^t\!\!,\hat{X}^t) &\geq \mathbb{E}_{P_{Y^{t-1}\!\!,U^{t-1}}}D_t(I(X^t;\hat{X}^t|y^{t-1}\!\!, u^{t-1})) \\
& \stackrel{(a)}{\geq} D_t(\mathbb{E}_{P_{Y^{t-1}\!\!,U^{t-1}}} I(X^t;\hat{X}^t|y^{t-1}\!\!, u^{t-1})) \\
&\stackrel{(b)}{\geq}  D_t(\gamma_t). 
\end{align*}
Recall that distortion-rate functions are in general convex and non-increasing
\cite[Lemma 10.4.1]{CoverThomas}. Thus, (a) follows from Jensen's inequality, and (b) follows from \eqref{eqgammat}.
Hence, under our privacy notion, \eqref{eqlossatt} ensures that the estimation error corresponding to any estimator $\hat{X}^t$ based on all information available in the cloud at time $t$ cannot be smaller than the distortion-rate function $D_t(\gamma_t)$.

{\upd
\subsection{Implication via Fano's inequality}
Suppose that $\mathcal{X}_t$ is a countable space for each $t=1, 2, ... , T$.
If we define $\epsilon_t\triangleq \text{Pr}(X^t\neq \hat{X}^t)$ under a joint distribution $P_{X^t,\hat{X}^t|y^{t-1},u^{t-1}}$, by Fano's inequality \cite{CoverThomas}, we have
\begin{align*}
&\epsilon_t\log(|\mathcal{X}^t|-1)+h(\epsilon_t) \\
&\geq H(X^t|\hat{X}^t,y^{t-1},u^{t-1}) \\
&=H(X^t|y^{t-1},u^{t-1})-I(X^t;\hat{X}^t|y^{t-1},u^{t-1})
\end{align*}
where $h(\epsilon)=-\epsilon\log \epsilon -(1-\epsilon)\log (1-\epsilon)$.
Taking the expectation with respect to $P_{Y^{t-1}, U^{t-1}}$,
\begin{align*}
&\mathbb{E}(\epsilon_t)\log (|\mathcal{X}^t|-1)+h(\mathbb{E}(\epsilon_t)) \\
&\geq \mathbb{E}(\epsilon_t)\log (|\mathcal{X}^t|-1)+\mathbb{E}(h(\epsilon_t)) \\
& \geq \mathbb{E} H(X^t|y^{t-1},u^{t-1})-\mathbb{E} I(X^t;\hat{X}^t|y^{t-1},u^{t-1})
\end{align*}
where we used Jensen's inequality in the first step.
Thus, from \eqref{eqgammat}, we obtain
\[
\mathbb{E}(\epsilon_t)\log (|\mathcal{X}^t|-1)+h(\mathbb{E}(\epsilon_t))
\geq H(X^t|Y^{t-1}, U^{t-1})-\gamma_t.
\]
This inequality clearly illustrates how the average probability of error is prevented from taking on small values by the randomness introduced, on average, when mapping $(Y^{t-1},U^{t-1})$ to $X^t$, as characterized by the conditional entropy $H(X^t|Y^{t-1},U^{t-1})$. The only way that the error probability can be allowed to become small is to counteract the growing uncertainty by increasing the value for $\gamma_t$, which illustrates how $\gamma_t$ captures loss of privacy by opening up for improved estimation of $X^t$ based on $(y^{t-1},u^{t- 1})$.
}

\section{LQG case}
\label{seclqg}

In this section, we consider a special case in which \eqref{minIstC} becomes a tractable optimization problem. Suppose the plant in Fig.~\ref{fig:out} is a fully observable linear dynamical system 
\[
X_{t+1}=A_tX_t+B_tU_t+W_t, \;\; Z_t=X_t
\]
where $W_t\sim\mathcal{N}(0,\Sigma^W_t)$ is a sequence of independent Gaussian random variables. 
{\upd We assume $\Sigma^W_t\succ 0$ for $t=1, ... , T$.}
Assume also that $c(\cdot, \cdot)$ in \eqref{minIstC} is a convex quadratic function, and that the problem \eqref{minIstC} can be written as
\begin{subequations}
\label{minIstCLQG}
\begin{align}
\min \;\; & I(X^T\rightarrow Y^T \| U^{T-1}) \label{minIstCLQGa} \\
\text{s.t.} \;\;\; & \sum_{t=1}^T \mathbb{E} (\|X_{t+1}\|_{Q_t}^2+\|U_t\|_{R_t}^2) \leq \delta. \label{minIstCLQGb}
\end{align}
\end{subequations}
The domain of optimization is \eqref{eqoptdomain}.\footnote{Problem \eqref{minIstCLQG} is identical to the problem considered in \cite{tanaka2015sdp}, except that in \cite{tanaka2015sdp}, an optimal solution is provided under the restriction that the stochastic kernels in \eqref{eqoptdomain} are Linear-Gaussian.} In what follows, we provide an optimal joint controller and output privacy filter policy that solves \eqref{minIstCLQG}.

\begin{figure}[t]
    \centering
    \includegraphics[width=0.9\columnwidth]{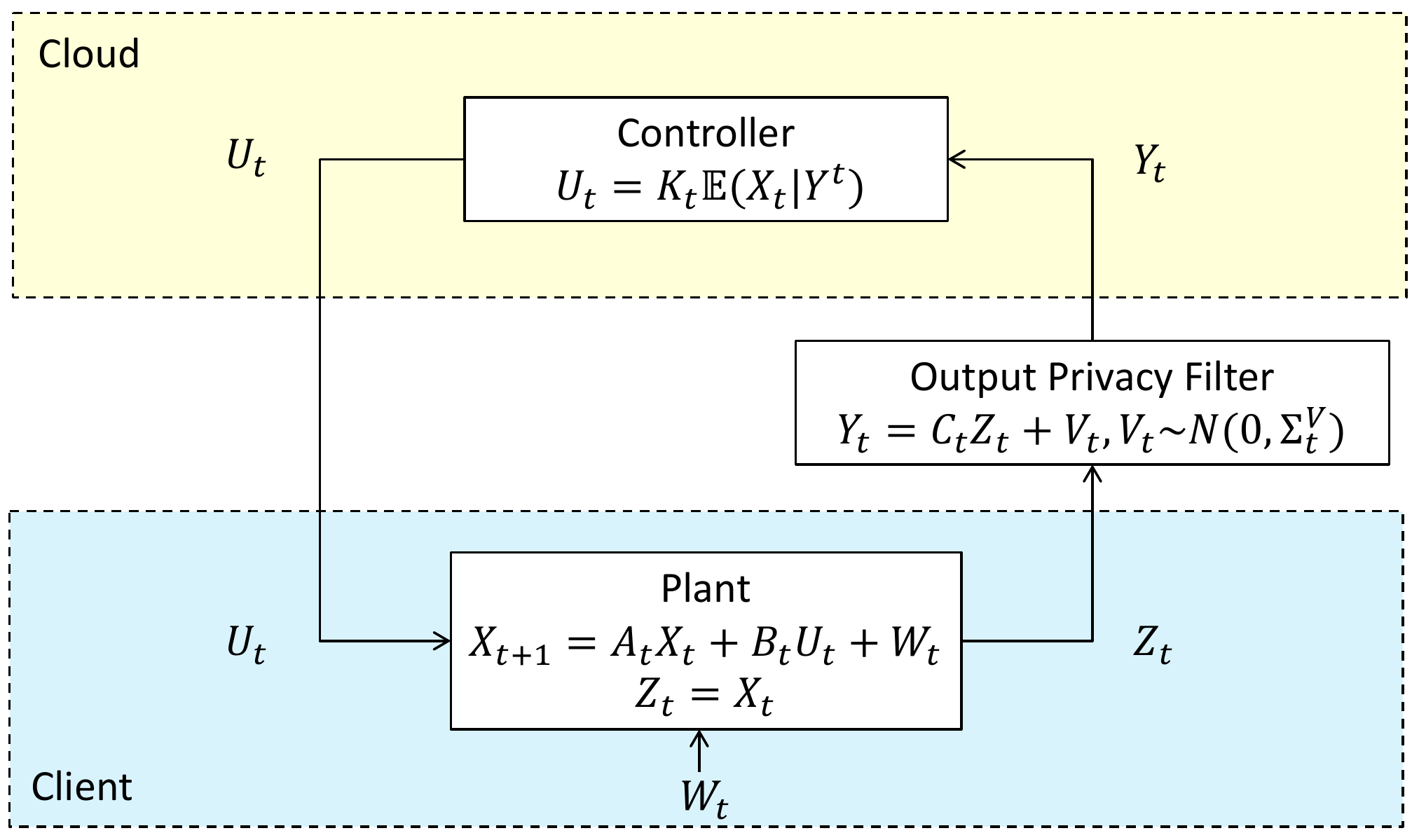}
    \caption{
    Structure of the optimal joint controller and output privacy filter for the cloud-based LQG control problem \eqref{minIstCLQG}. Although the output privacy filter is allowed to utilize public random variable $U^{t-1}$ (as shown in Fig.~\ref{fig:out}), it turns out that this information need not be used.}
    \label{fig:outlqg}
\end{figure}

{\upd
First, in view of Lemma~\ref{lemfbineq}, notice that the minimum privacy leakage characterized by \eqref{minIstCLQG} is lower bounded by the optimal value of
\begin{subequations}
\label{mindi}
\begin{align}
\min \;\; & I(X^T\rightarrow U^T ) \\
\text{s.t.} \;\;\; & \sum_{t=1}^T \mathbb{E} (\|X_{t+1}\|_{Q_t}^2+\|U_t\|_{R_t}^2) \leq \delta
\end{align}
\end{subequations}
where again the domain of optimization is $\mathcal{D}$ given by \eqref{eqoptdomain}.

A related optimization problem to \eqref{mindi} is already considered in \cite{tanaka2015lqg}, where the only difference is that the optimization domain considered there is $\mathcal{D}'\triangleq \{P_{U_t|X^t,U^{t-1}}\}_{t=1}^T$. Notice that $\mathcal{D} \subset \mathcal{D}'$ since every element in $\mathcal{D}$ can be, by compositions of stochastic kernels, mapped to an element of $\mathcal{D}'$.
In \cite{tanaka2015lqg}, it is shown that the optimal solution $P_{U_t|X^t,U^{t-1}}$ in $\mathcal{D}'$ can be realized by the interconnection of a linear sensor, Kalman filter, and a controller as shown in Fig.~\ref{fig:threestage}. Matrix parameters such as $C_t$, $\Sigma_t^V$, $L_t$, $K_t$ for $t=1, ... , T$ must be optimally tuned, which can be achieved by an LMI algorithm (Algorithm~\ref{algsdp}). Notice that in the second step of Algorithm~\ref{algsdp}, a convex optimization problem must be solved. Since it is in the form of the determinant maximization problem, a standard semidefinite programming solver can be used.

\begin{figure}[t]
    \centering
    \includegraphics[width=0.9\columnwidth]{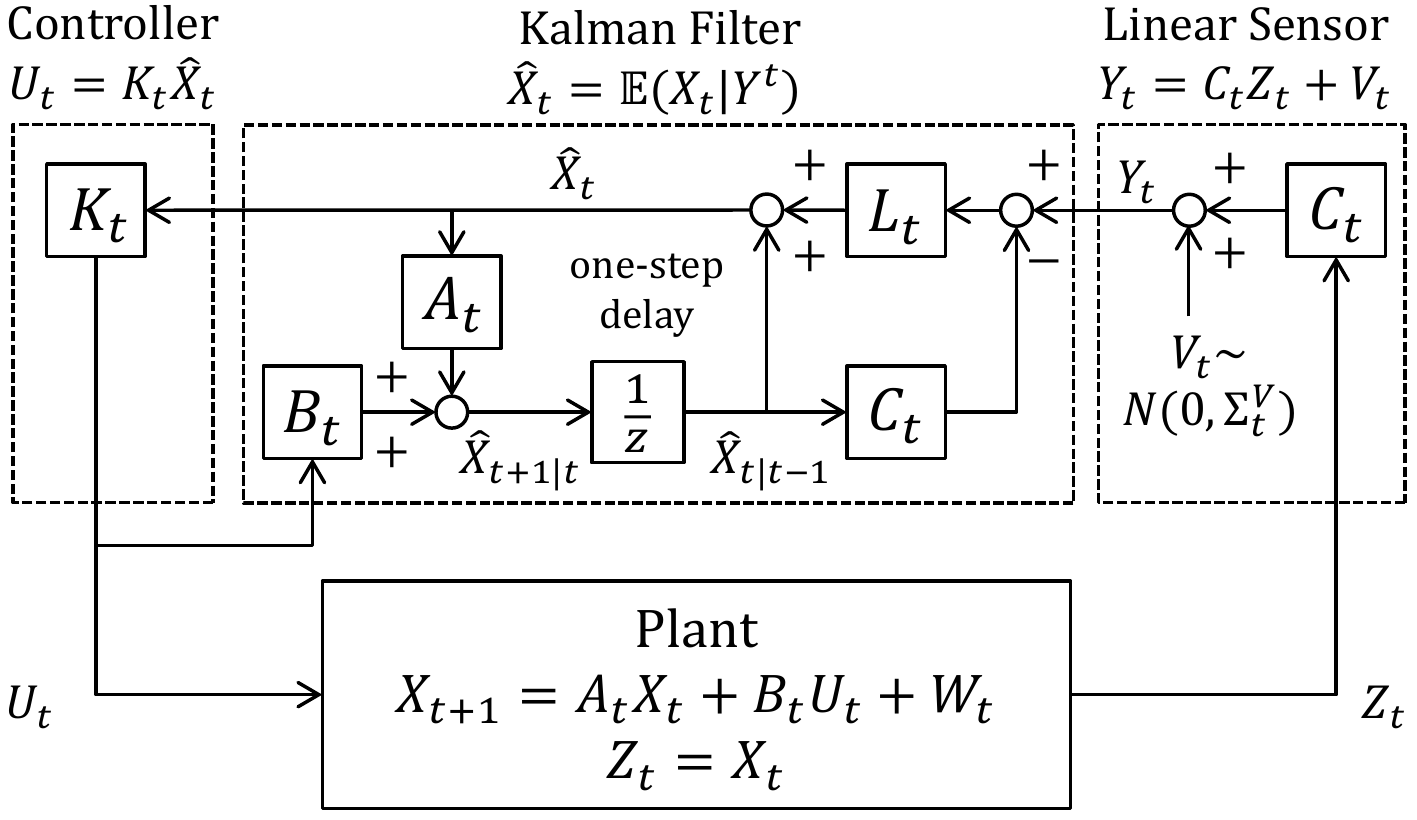}
    \caption{Structure of optimal policy for problem \eqref{mindi}. }
    \label{fig:threestage}
\end{figure}

For the purpose of our present study, notice that the linear sensor mechanism in Fig.~\ref{fig:threestage} can be interpreted as an output privacy filter, while the Kalman filter and the controller can be implemented by the cloud operator (Fig.~\ref{fig:outlqg}). In other words, a realization of $P_{U_t|X^t,U^{t-1}}$ shown in Fig.~\ref{fig:threestage} can be viewed as an element of $\mathcal{D}$. 
}
Now we claim the following.

\begin{lem}
\label{proplqg}
The policy shown in Fig.~\ref{fig:threestage}, which is an optimal solution to \eqref{mindi}, is also an optimal solution to  \eqref{minIstCLQG}.
\end{lem}
{\upd
\begin{pf}
Since $I(X^T\rightarrow U^T)\leq I(X^T\rightarrow Y^T\|U^{T-1})$ by Lemma~\ref{lemfbineq}, it is sufficient to show that 
\begin{equation}
\label{eqdixyu}
I(X^T\rightarrow U^T)= I(X^T\rightarrow Y^T\|U^{T-1})
\end{equation}
holds for the optimal solution \eqref{mindi} shown in Fig.~\ref{fig:threestage}. 
For the ease of presentation, we shown \eqref{eqdixyu} under the assumption that $K_t$, $t=1, ... , T$ are invertible. Although such an assumption is unnecessary, the proof in full generality is lengthy and must be differed to \cite[equation (22)]{tanaka2015lqg}.

Due to the invertibility of $K_t$, $t=1, ... , T$, random variables $U^t$ and $\hat{X}^t$ appearing in Fig.~\ref{fig:threestage} are related by an invertible linear map. Since $C_t$, $t=1, ... , T$ are full row rank matrices by construction, the Kalman filter
\[
\hat{X}_t=(I-L_tC_t)(A_{t-1}\hat{X}_{t-1}+B_{t-1}U_{t-1})+L_tY_t
\]
has a causal inverse
\[
Y_t=L_t^\dagger \hat{X}_t+L_t^\dagger (L_tC_t-I)(A_{t-1}\hat{X}_{t-1}+B_{t-1}U_{t-1})
\]
where $L_t^\dagger\triangleq (L_t^\top L_t)^{-1}L_t^\top$. Thus, $\hat{X}^t$ and $Y^t$ are related by an invertible linear map. Therefore, $U^t$ and $Y^t$ are related by an invertible linear map. (They contain statistically equivalent information.) In particular, this implies that conditional differential entropies $h(X^t|Y^{t-1},U^{t-1})$ and $h(X^t|Y^t,U^{t-1})$ are equal to $h(X^t|U^{t-1})$ and $h(X^t|U^t)$, respectively. Therefore,
\begin{align*}
&I(X^T\rightarrow Y^T\|U^{T-1}) \\
&=\sum\nolimits_{t=1}^T I(X^t;Y_t|Y^{t-1},U^{t-1}) \\
&=\sum\nolimits_{t=1}^T \left[h(X^t|Y^{t-1},U^{t-1})-h(X^t|Y^t,U^{t-1})\right] \\
&=\sum\nolimits_{t=1}^T \left[h(X^t|U^{t-1})-h(X^t|U^t)\right] \\
&=\sum\nolimits_{t=1}^T I(X^t;U_t|U^{t-1}) \\
&=I(X^T\rightarrow U^T). \hspace{30ex}\text{ \QED}
\end{align*}
\end{pf}
}
In summary, we have an explicit form of the joint control and output privacy filter policy solving \eqref{minIstCLQG}.
\begin{prop}
\label{propoptjoint}
An optimal joint controller and output privacy filter characterized by an optimal solution to \eqref{minIstCLQG} is in the form shown in Fig.~\ref{fig:outlqg}.
An optimal choice of matrices $C_t$, $\Sigma_t^V$, $L_t$ (Kalman gains) and $K_t$ (feedback control gains) are obtained by Algorithm~\ref{algsdp}.
Moreover, the optimal value of \eqref{minIstCLQG} is equal to the optimal value of the determinant maximization problem in Algorithm~\ref{algsdp}.
\end{prop}

\begin{algorithm}
  \caption{Joint controller and privacy filter design for cloud-based LQG control}
  {\upd
\begin{enumerate}[label=$\bullet$, itemsep = 2mm, topsep = 2mm, leftmargin = 3mm]
\item[1.]  Determine  feedback control gains $K_t$ via the backward Riccati recursion:
\begin{align*}
S_t&=\begin{cases} Q_t & \text{ if } t=T \\ Q_t+\Phi_{t+1} & \text{ if } t=1,\cdots, T-1 \end{cases} \\
\Phi_t&=A_t^\top (S_t-S_t B_t (B_t^\top S_t B_t + R_t)^{-1} B_t^\top S_t) A_t \\
K_t&= -(B_t^\top S_t B_t + R_t)^{-1} B_t^\top S_t A_t \label{backricK}\\
\end{align*}
\item[2.] Solve a determinant maximization problem with respect to $P_{t|t} \succ  0, \Pi_t\succ  0$, $t=1, ... , T$ subject to LMI constraints:
\begin{align}
\min & \quad \frac{1}{2} \sum\nolimits_{t=1}^T \log\det \Pi_t^{-1} + c_1 \nonumber \\
\text{s.t.} \;& \quad \sum\nolimits_{t=1}^T \text{Tr}(\Theta_t P_{t|t}) + c_2 \leq D,  \nonumber \\
& \quad  P_{1|1}\preceq P_{1|0}, P_{T|T}=\Pi_T,  \nonumber \\
& \quad  P_{t+1|t+1}\preceq A_t P_{t|t}A_t^\top +\Sigma_t^W, \;\; t=1, ... , T-1 \nonumber \\
&\hspace{1ex} \left[\! \begin{array}{cc}P_{t|t}\!-\!\Pi_t\!\! &\!\! P_{t|t}A_t^\top \nonumber  \\
A_tP_{t|t} \!&\! A_t P_{t|t}A_t^\top\!+\!\Sigma_t^W  \end{array}\!\right]\! \succeq\! 0, \;\;  t=1, ... , T-1 \nonumber
\end{align}
where $\Theta_t= K_t^\top (B_t^\top S_t B_t + R_t) K_t$, $t=1, ... , T$ and
\begin{align*}
c_1&=\tfrac{1}{2}\log\det{P_{1|0}}+\tfrac{1}{2}\sum\nolimits_{t=1}^{T-1} \log \det \Sigma_t^W \\
c_2&=\text{Tr}(N_1P_{1|0})+\sum\nolimits_{t=1}^T \text{Tr}(\Sigma_t^W S_t).
\end{align*} 
\item[3.] For each $t=1, ... , T$, choose (e.g., by the singular value decomposition) a full row rank matrix $C_t$ and a positive definite matrix $\Sigma^V_t$ such that  
\[
C_t^\top {\Sigma^V_t}^{-1}C_t=P_{t|t}^{-1}-(A_{t-1}P_{t-1|t-1}A_{t-1}^\top + \Sigma_{t-1}^W)^{-1}.
\]
\item[4.]  Determine the Kalman gains by
\[
L_t=P_{t|t-1}C_t^\top (C_t P_{t|t-1}C_t^\top + \Sigma^V_t)^{-1}
\]
where $P_{t+1|t}=A_tP_{t|t}A_t^\top+\Sigma_t^W$.
\end{enumerate}
}
  \label{algsdp}
\end{algorithm}

Notice that the privacy filter shown in Fig.~\ref{fig:outlqg} is similar to privacy protecting mechanisms considered in various other contexts (e.g., \cite{dwork2008differential}) in that it is adding noise $V_t$ before disclosing data. Proposition~\ref{propoptjoint} shows that the optimal noise distribution is Gaussian in the LQG case \eqref{minIstCLQG}.

\section{Numerical example}

\begin{figure}[t]
    \centering
    \includegraphics[width=\columnwidth]{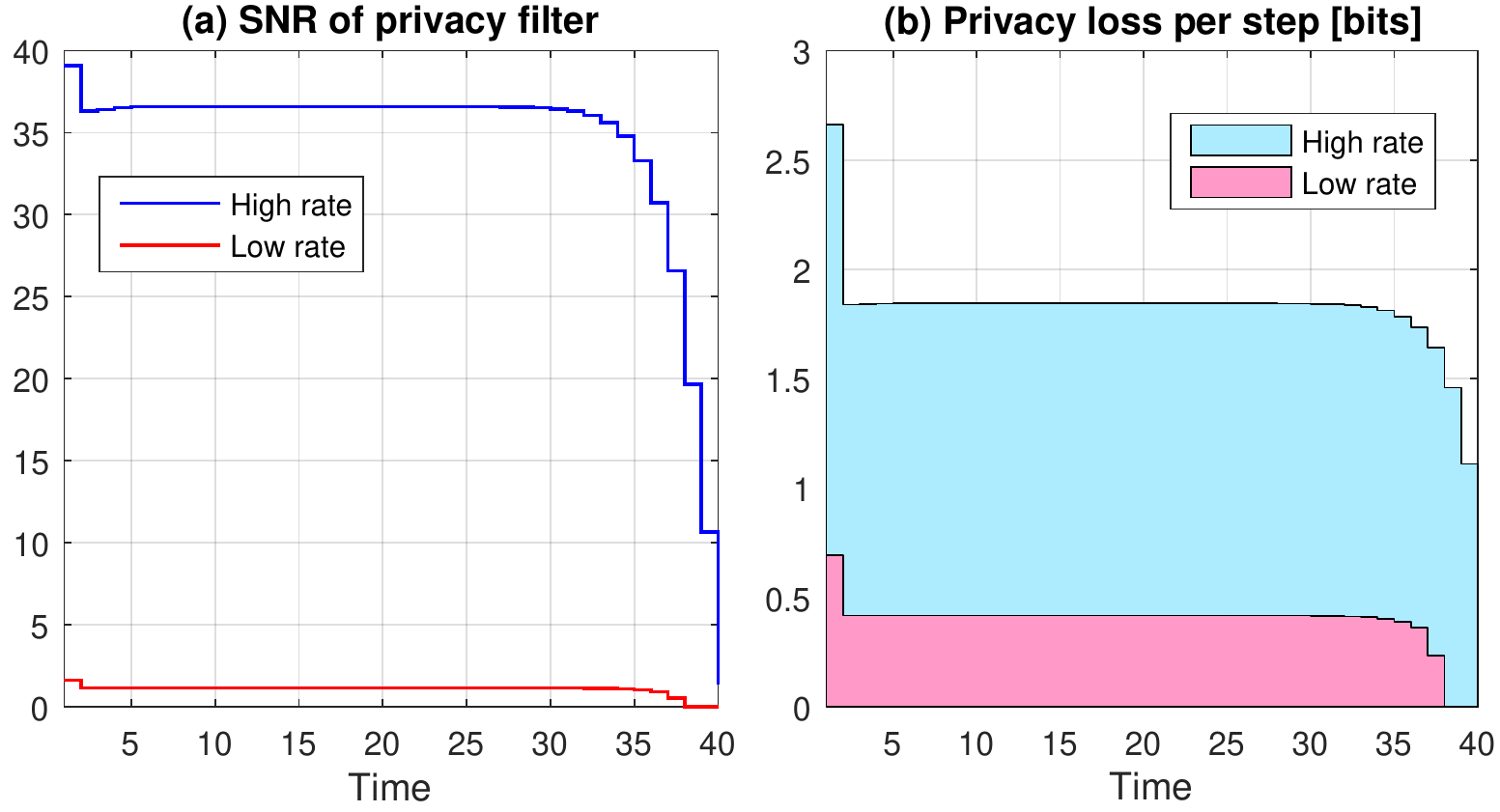}
    \caption{SNR of the optimal privacy filter and the privacy loss.}
    \label{fig:snr}
\end{figure}

\begin{figure}[t]
    \centering
    \includegraphics[width=\columnwidth]{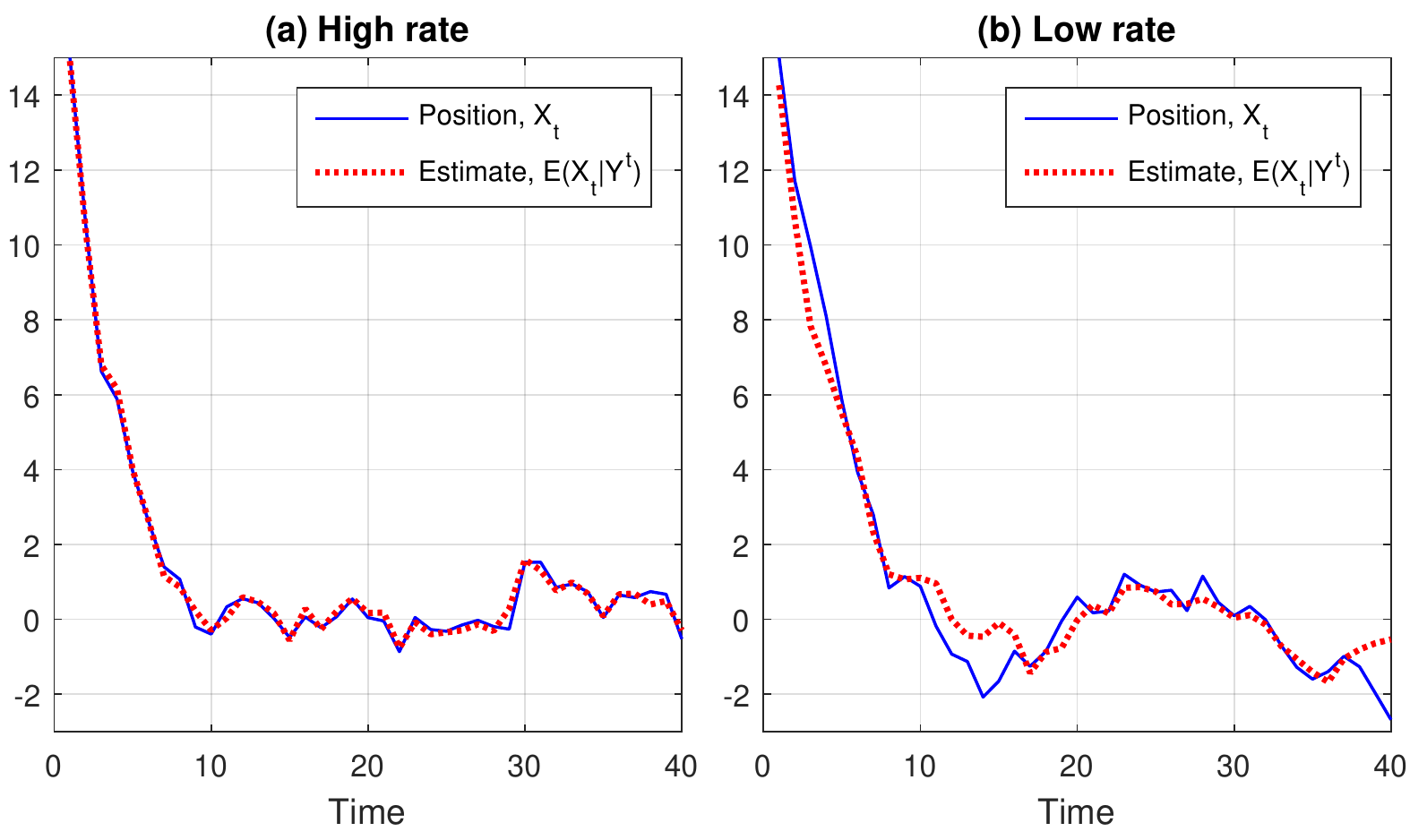}
    \caption{Simulated trajectory of $X_t$ and its estimate $\mathbb{E}(X_t|Y^t)$.}
    \label{fig:x_est}
\end{figure}

In this section, we consider a simple scalar system
\[ X_{t+1}=X_t+U_t+W_t, \;\; t=1, ... , T\]
with a process noise $W_t\sim\mathcal{N}(0, 0.3)$. This example is motivated by a cloud-based navigation service, where the state variable $X_t$ is interpreted as the position of the client at time $t$, whereas $U_t$ is the navigation signal provided by the cloud.
Assuming the initial position is {\upd $X_1=15$}, the cloud-based controller navigates the client to the origin withing {\upd $T=40$} steps using  the output of a privacy filter
\[ Y_t=C_tX_t+V_t, \;\; V_t\sim\mathcal{N}(0, \Sigma_t^V).\]
The optimal privacy filter is different depending on the choice of $\delta$ in \eqref{minIstCLQG}.
We consider two scenarios in which control requirements are stringent {\upd ($\delta=24.4$)} and mild {\upd ($\delta=31.4$)}. In both cases, we use the same control cost function with {\upd$Q_t=1, R_t=10$} for all $t$.
The former case requires higher data rate (measured in directed information).
In each scenario, we compute the optimal sequence $\{C_t, \Sigma_t^V\}_{t=1, ... , 40}$ by solving \eqref{minIstCLQG} using semidefinite programming.
The sequence of signal-to-noise ratios $\text{SNR}_t=C_t^2/\Sigma_t^V$ is plotted in Fig.~\ref{fig:snr} (a). The total loss of privacy in the high rate scenario is {\upd $71.5$} [bits]  (the area of blue region in Fig.~\ref{fig:snr} (b)) while it is {\upd $15.4$} [bits]  in the low rate case (the area of red region in Fig.~\ref{fig:snr} (b)).
Fig.~\ref{fig:x_est} shows the closed-loop performance in each scenario. In the high rate case, the cloud estimates the position of the client accurately, and consequently the control performance is better. In the low rate case, the estimate is not accurate (privacy is better protected) so the control performance is poor.

\section{Summary and future work}
\label{secsummary}
In this paper, we showed that Kramer's causally conditioned directed
information arises as the unique candidate (up to positive
multiplicative factor) for a measure of privacy loss in cloud-based
control if Postulates~\ref{post1}-\ref{post4} (including \emph{the
  data-processing axiom} \cite{jiao2015justification}) are to be
satisfied. Our result is a first step towards a complete axiomatic
privacy theory in networked control, which we see as a promising
approach that provides a convenient interface between theory and
practice.  There are numerous further opportunities along the same line
of research. Notice that the set of postulates we have selected in
this paper is not the only possible characterization of privacy. In
fact, in \cite{du2012privacy}, the ``maximum'' type of privacy leakage
function
\begin{align*}
&L(\ell,P_{X,Y})= \\
& \;\; \inf_{Q_X} \mathbb{E}_{P_{X}} [\ell(X, Q_X)]-\min_{y\in \mathcal{Y}}\inf_{Q_{X|Y}} \mathbb{E}_{P_{X|Y}} [\ell(X, Q_{X|Y})].
\end{align*}
is considered in parallel with the ``average'' type of privacy leakage
function we assumed in Postulate~\ref{post1}. It is of great interest
whether there exists a valid notion of privacy satisfying the
corresponding new set of postulates.


\bibliographystyle{plain}        
\bibliography{ref}           



\end{document}